\documentclass[review]{elsarticle}
\usepackage[bookmarks=true,colorlinks=true,urlcolor=blue,citecolor=red,linkcolor=blue,linktocpage,pdfpagelabels,bookmarksnumbered,bookmarksopen]{hyperref}
\usepackage{lineno,hyperref}

\modulolinenumbers[5]

\journal{arXiv}
\usepackage{amscd}
\usepackage{mathrsfs}
\usepackage{amsfonts}
\usepackage{graphicx}
\usepackage{amsmath,amscd,stmaryrd,amssymb,array, amsfonts,amsmath,amssymb}
\usepackage{enumitem}
\usepackage{color}
\usepackage{verbatim}
\usepackage{extarrows}
\usepackage{amsmath}
\usepackage{cases}

\numberwithin{figure}{section}
 \numberwithin{equation}{section}

\newtheorem{theorem}{Theorem}[section]
\newtheorem{proposition}[theorem]{Proposition}
\newtheorem{definition}[theorem]{Definition}

\newtheorem{lemma}[theorem]{Lemma}
\newtheorem{remark}[theorem]{Remark}

\newcommand{\bbf}{{f}}

\newcommand{\bh}{{h}}

\newcommand{\bV}{{ V}}

\newcommand{\cC}{{\mathcal C}}

\newcommand{\cL}{{\mathcal L}}

\newcommand{\sX}{{\mathscr X}}

\def\be{\begin{equation}}
\def\ee{\end{equation}}
\def\bes{\begin{equation*}}
\def\ees{\end{equation*}}
\def\bsp{\begin{split}}
\def\esp{\end{split}}

\def\ba{\begin{array}}
\def\ea{\end{array}}
\def\benu{\begin{enumerate}}
\def\eenu{\end{enumerate}}
\def\bt{\begin{theorem}}
\def\et{\end{theorem}}
\def\bp{\begin{proposition}}
\def\ep{\end{proposition}}
\def\bl{\begin{lemma}}
\def\el{\end{lemma}}
\def\br{\begin{remark}}
\def\er{\end{remark}}
\def\bd{\begin{definition}}
\def\ed{\end{definition}}

\def\b{\beta}
\def\De{\Delta}

\def\pa{\partial}

\def\lam{\lambda}
\def\Lam{\Lambda}

\def\gam{\gamma}

\def\W{\Omega}

\def\.{\cdot}

\def\R{\mathbb{R}}

\def\ol{\overline}

\def\ra{\rightarrow}

\def\~{\tilde}
\def\8{\infty}
\def\X{\times}

\def\mb{\mbox}

\def\llg{\left\langle}
\def\rrg{\right\rangle}

\def\Hs{\hspace{1cm}}\def\hs{\hspace{0.5cm}}
\def\Vs{\vskip8pt}\def\vs{\vskip4pt}

\def\({\left(}\def\){\right)}

\begin{document}

\begin{frontmatter}

\title{Global Existence for Reaction-diffusion Equations with State-Dependent Delay and Fast-growing Nonlinearities\footnote{This work was supported by Shandong Provincial Natural Science Foundation [No. ZR2024QA228], National Natural Science Foundation of China [12171352]. }
}


\author[mymainaddress]{Ruijing Wang }
\ead{wrj\_math@tju.edu.cn}
\address[mymainaddress]{School of Science, Qingdao University of Technology, Qingdao 266520,  China}


\begin{abstract}
This work aims to study the initial-boundary value problem of the reaction-diffusion equation with state-dependent delay
$\pa_{t}u-\Delta u=f(u)+g(u,u(t-\tau(t,u_t)))+h(t,x)$ in a bounded domain.
We establish the global existence of the problem under suitable dissipative-type structural conditions, allowing both nonlinear terms $f$ and $g$ to have arbitrary polynomial growth rates. Another highlight in this work is that, we significantly relax  the continuity assumptions imposed  on the delay functions.
\end{abstract}


\end{frontmatter}


\section{Introduction}\label{s:1}

We revisit  the global existence of the following  initial-boundary value problem
\be\label{e:1.1}\left\{\ba{ll}
\pa_{t}u-\Delta u=f(u)+g(u,\,u(t-\tau(t,u_t)))+h(t,x),\\[1ex]
u|_{\pa\W}=0,\hs u|_{\overline{\W}\times [-r,0]}=\phi\ea\right.
\ee
 with state-dependent delays. Here $\W\subset \R^d$ is a bounded domain with sufficiently smooth boundary $\pa\W$,  $f$ and $g$ are assumed to be  $C^1$ functions,  $r\geq 0$ is a fixed number, and $\tau$ denotes the time delay which may depend on the history $u_t$ of the state $u(\cdot)$ at the time $t$ given by
\be\label{e:1.5}
u_t(\theta)=u(t+\theta), \Hs \theta \in [-r,0].
\ee

Functional differential equations with state-dependent delays (SDDs) were initially introduced in a much earlier work \cite{Poi} by Poisson.
But they were not taken seriously until Driver \cite{Dri} established  a more  intuitive  model to describe a  two-body problem in electrodynamics.
Since then, an increasing number of differential equations  with SDDs have been proposed to deal with various problems in different areas such as stage-structured population growth, automatic position control, the dynamics of white blood cells, and synaptic delayed adaptive processes of neurons; see e.g. \cite{Aie,Cal, Har,LXD} and references therein. This aroused much interest in the study  of such systems over the past decades; see e.g. \cite{Bar,Chu,Har,Hu-Wu,Kri1,Kri, KW2,Lv,MN4,QW,Rez3,Rez2,Wal1,Wal2,Wu}, to name a few.

The aim of this paper is to study the global existence, regularity, and dissipativity of the PDE system \eqref{e:1.1}.
These questions  have already been addressed in Rezounenko \cite{Rez3},  Krisztin and Rezounenko \cite{Kri}, Hern\'{a}ndez, Fernandes and  Wu \cite{Her1,Her4} etc. in certain abstract functional framework, where one can find some generally theory that can be successfully applied to parabolic PDEs with SDDs and sublinear or subcritical nonlinearities. Here we are interested in the case the nonlinear terms  $f$ and $g$ in \eqref{e:1.1} are permitted to  have arbitrary polynomial growth rates:
 \benu
  \item[{\bf(H0)}] There exist positive constants $p,\b,a_0$ and $b_0$ such that
  $$
|f(u)|\leq a_0(|u|^p +1),\hs|g(u,v)|\leq b_0(|u|^\beta +|v|^\beta +1),\Hs u,v\in \mathbb{R}.
$$
\eenu
Since there is no restriction on the exponents $p$ and $\b$, system \eqref{e:1.1} can be of supercritical.  In such a case the system  may fall outside the scope of the existing theoretical results in the literature mentioned above. Another obstacle to study the global well-posedness for such a  system  is that, the fast-growing nonlinearities may cause finite-time blowup phenomenon and destroy the uniqueness of solutions. In fact, even if for a reaction-diffusion system without delay, the global well-posedness and dynamics is a challenging problem when the system involves supercritical nonlinearity; see e.g. the recent paper  \cite{Kostianko} by Kostianko et al.
To overcome the difficulty mentioned above, we will work in a suitable framework of  solutions in certain weaker sense.

For this purpose, the first step is to find appropriate structure conditions on the nonlinearities  and establish necessary global estimates for weak solutions. This is in general not a easy task. A widely used class of structure conditions on the nonlinearities is of dissipative type, which have been extensively applied to both delay and non-delay equations in the literature, leading to nice results in many aspects. Here, inspired by a recent work \cite{Li2}, we impose on \eqref{e:1.1}  the following disspative-type hypothesis:
\benu
  \item[{\bf(H1)}] There is  $\b_0>\b$ and $\Lam,N> 0$ such that
$$
\bbf(s)s\leq -\Lam |s|^{{\b_0}+1}+N,\Hs s\in\R.
$$
\eenu

For notational simplicity, given a Banach space $X$, we write
\be\label{e:1.6}
\mathcal{C}_{X}=C([-r,0];X),\hs \mathcal{L}^{\infty}_{X}=L^{\infty}(-r,0;X).
\ee
 Here $r$ is the upper bound of the delay function $\tau$ appearing in \eqref{e:1.1}.
To guarantee well-posedness of state-dependent delay differential equations like \eqref{e:1.1}, in the literature one usually imposes on the delay functions  a set of continuity and structure conditions, which seem to be a little complicated and stronger; see e.g.  \cite{Car2,Her4,Her5}. In this work  we only  assume   that
\benu
  \item[{\bf(H2)}] $\tau\in C\(\R\X \cC_{L^2(\W)};[0,r]\)$ with $\tau(t,\phi)$ being locally Lipschitz continuous in $\phi\in\cC_{L^2(\W)}$ in a uniform manner with respect to $t$ in any compact interval.
\eenu
Under the above hypotheses (H0)-(H2), we show that there exists a critical exponent $q_{c}>0$ such that if $q>q_{c}$, then for every initial data $$\phi\in \cC_{H^1_0(\W)}\cap \cL^\8_{L^q(\W)}:=\sX_1^q,$$
  system \eqref{e:1.1} has a global weak solution $u$ satisfying appropriate  exponential decay estimates in $L^q(\W)$ and $H^1(\W)$; see Theorem \ref{t:3.1} for details. Furthermore, the system has an instantaneous $L^{\8}$-smoothing property.

This paper is organized as follows. In Section \ref{s:2} we recall several functional spaces and give definitions of weak solutions for \eqref{e:1.1}.
Section \ref{s:3} is devoted to the existence and global $L^{q}$, $H^{1}$ and $L^{\infty}$-estimates for the weak solution $u$ of \eqref{e:1.1}.

\section{Preliminaries}\label{s:2}
\subsection{Functional  spaces}\label{s:2.2}
\noindent $\bullet$ {\bf Sobolev spaces and  their associated functional spaces.}
The spaces $W^{k,p}(\W)$ and $W^{k,p}_0 (\W)$ are standard Sobolev spaces. Denote by $|\cdot|_p$ the usual norm on the Lebesgue space $L^p (\W)$, and let $(\cdot,\cdot)$ be the scalar product on $H=L^2(\Omega)$.
Let $H^k(\W)=W^{k,2}(\W)$, and $H^k_0(\W)=W^{k,2}_0(\W)$. Set the Hilbert spaces
$$
V_1=H^1_0(\W),\hs V_2=H^2(\W)\cap H^1_0(\W)
$$
with  the following scalar products respectively:
$$
\llg u,v \rrg_1=(\nabla u,\nabla v) \,\, (u,v\in V_1),\hs \llg u,v \rrg_2=(  \Delta u,  \Delta v) \,\, (u,v\in V_2).
$$
Denote by $\|\cdot\|_{V_k}$ ($k=1,2$) the norm induced by the scalar product $\llg \cdot,\cdot \rrg_k$. It is trivial to check that $\|\cdot\|_{V_k}$ is a norm on $V_k$ equivalent to the usual $H^k(\W)$-norm.

Let $r$  be the upper bound of the delay function $\tau$ appearing in \eqref{e:1.5}. Given a Banach space $X$, let $\cC_X$ and $\cL_X^\8$ be the Banach spaces defined in \eqref{e:1.6},
which  are equipped respectively with the norms
$$
\|\cdot\|_{\mathcal{C}_{X}}:=\|\cdot\|_{C([-r,0];X)},\hs \|\cdot\|_{\mathcal{L}^{\infty}_{X}}:=\|\cdot\|_{L^{\infty}(-r,0;X)}.
$$
In the case $r=0$,  we also assign
$$\mathcal{C}_{X}=\mathcal{L}^{\infty}_{X}=X.$$

\subsection{Definition of a weak solution}\label{s:2.1}
Let $f,g,h$ be the functions in \eqref{e:1.1}. We may write  $h(t)=h(t,\.)$, and let
$$
f(u)+g(u,u(t-\tau(t,u_t)))+h(t):=F(t,u,u_t).
$$
Then system \eqref{e:1.1}  can be put into an abstract one:
\begin{numcases}{}
       u'+Au=F(t,u,u_t),\hs t>0,\label{E2.1}\\
      u(s)=\phi(s),\hs s\in [-r,0].\label{E2.2}
\end{numcases}

Let $p$ and $\b$ be the exponents in (H0), and set
$$\hat{p}=\max\{p,\b\}.$$
Let $J\subset \mathbb{R}$ be an interval.
\begin{definition}\label{d:1.1}
 A function $u=u(t)$ taking values in $V_1$ is a weak solution of \eqref{E2.1} on $J$, if for any compact interval $[a,b]\subset J$,
\be\label{e:2.5}
\begin{split}
u\in C([a-r,b];H)\,+\, &L^2(a-r,b;V_1)\,+ L^{2\hat{p}}(a-r,b;L^{2 \hat{p}}(\W)),
\\[.5ex]
&u'\in L^2(a,b;V'_1),
\end{split}
\ee
and it holds in the distribution sense on $J$ that
$$
\llg u',w\rrg+\llg u,w\rrg_1=\(F(t,u,u_t),\,w\),\Hs \forall w\in V_1,
$$
where $\llg \cdot,\cdot \rrg$ stands for the paring between $V_1$ and $V'_1$.
\end{definition}

\begin{remark}\label{r:1.1}
By definition, it is easy to verify that if $u$ and $v$ are weak solutions of \eqref{E2.1} on the intervals $[t_0, t_1]$ and $[t_1, t_2]$ respectively with
$$u(t)=v(t),\Hs t \in [t_1-r, t_1],$$ then the union $z$ of $u$ and $v$ is a weak solution of \eqref{E2.1} on $[t_0, t_2]$, where
$$
z(t)=u(t)\,(t\leq t_1),\hs z(t)=v(t)\,(t>t_1).
$$
\end{remark}

\br The requirements on a weak solution $u$ of \eqref{E2.1} on $J$ and the growth restrictions on the nonlinear terms in (H0) imply that $F(t,u,u_t)$ belong to $L^2(a,b;H)$ for any compact interval $[a,b]\subset J$. Further by standard theory on linear evolution  equations (see e.g. \cite [Chap. 2, Section 3.2]{Tem}), it can be easily deduced that $u\in C(J,H)$.\er




\begin{definition}\label{d:1.2} Given $\phi\in\cC_H$, a weak solution $u$ of \eqref{E2.1} on $(0,T)$ is called a weak solution of system \eqref{e:1.1} (or instead, \eqref{E2.1}-\eqref{E2.2}) on $[-r,T]$, if $u\in C([-r,T);H)$ and fulfills the initial condition \eqref{E2.2}.
\end{definition}

In what follows,  we will simply call a \emph{weak solution} of \eqref{e:1.1} as a \emph{solution}, unless otherwise stated.

\section{Global Existence and Estimates of Weak Solutions}\label{s:3}
In this section we first derive some global estimates for the Galerkin approximations of the solutions of \eqref{e:1.1}. Then we  employ some compactness argument  to establish the global existence result. We also address the instantaneous $L^{\8}$-smoothing property for \eqref{e:1.1} and give a global $L^\8$-estimate.

Let $p$, $\b$ and $\b_0$ be the exponents  in (H0) and (H1). Set
\be\label{e:1.3}
q_{c}=\mb{max}\{2p,2\b,p_0\}, \hs  \mb{where }\, p_0=\frac{(\b_0-1)}{\b_0-\b} \b.
\ee
For notational convenience, we also write
\be\label{e:1.7}
\mathscr{X}^q_i=\mathcal{C}_{V_i}\cap \mathcal{L}^{\infty}_{L^q(\W)},\Hs i=1,2.
\ee

\subsection{Global $L^{q}$ and $H^1$-estimates of Galerkin approximations}\vs
Let $A=-\De$ be the Dirichlet Laplacian with eigenvalues
$$0<\mu_1\leq \mu_2\leq\cdots\leq\mu_j\leq\cdots.$$
 Denote $w_j$ ($j=1,2,\cdots$)  the corresponding normalized eigenfunctions which  form an orthogonal basis of $H:=L^2(\W)$.

Given $\phi \in \sX_1^q$, by \cite[Theorem 7.1]{Li2} one can pick a sequence of smooth functions $\phi_k=\sum^k_{j=1}c_{kj}(t)w_j$\,($k=1,2,\cdots$) such that
\be\label{e:2.1}
\phi_k\rightarrow\phi \,\,( \mb{in}\,\, \mathcal{C}_{V_1}),\hs \|\phi_k\|_{\mathcal{L}^{\8}_{L^q (\W)}}\leq 8 \|\phi\|_{\mathcal{L}^{\8}_{L^q (\W)}}\,\, (k\in \mathbb{N}).
\ee
Similarly we can take a sequence of functions $\{h_k\}$ with $h_k(t)=\sum^k_{j=1}d_{kj}(t)w_j$ such that $h_k \rightarrow h$ in appropriate functional  spaces. For each $k$, let
\be\label{e:2.2}
u_k(t)=\sum^k_{j=1}a_{kj}(t)w_j
\ee
be a Galerkin  approximation solution of \eqref{e:1.1},
\be\label{e:2.3}
\left\{
\begin{split}&\(\frac{d u_k}{dt}-\De u_k,w_j\)=\( F_k(t,u_k,(u_k)_t),\,w_j\),\hs 1\leq j\leq k,\\
 & u_k(s)=\phi_k(s),\hs s\in[-r,0],
\end{split}
\right.
\ee
where
\bes\label{e:x2}
F_k(t,u_k,(u_k)_t):=f( u_k)+g(u_k, u_k(t-\tau(t,(u_k)_t)))+ h_k.
\ees
By the basic theory of  ODEs, we know that \eqref{e:2.3} has a unique solution $u_k(t)=u_k(t;\phi)$. Furthermore,
$u_k$ is sufficiently regular so that all the calculations below can be performed rigorously on  $u_k$.




\bl\label{t:1.1}
Assume  {\em (H0)-(H2)}. Let $q_{c}< q< \8$. 
 Then there exist $\lam_q,R_q,\lam_1,\rho_1>0$ and $M_i$ ($i=0,1$) ($M_i's$ are  independent of $q$) such that for all $\phi\in \mathscr{X}_1^q$,
\be\label{e:3.1}
|u_k(t;\phi_k)|_{q}^{q}\leq M_0\|\phi_k\|_{\cL_{L^q(\W)}^\8}^{q}e^{-\lam_q t}+R_q,\Hs \,t\geq 0.
\ee
and
\be\label{e:3.6}
\| u_k(t;\phi_k)\|_{V_1}^{2} \leq  \|\phi_k\|_{\cC_{V_1}}^2e^{-\mu_1 t}+M_1\|\phi_k\|_{\cL_{L^q(\W)}^\8}^{q}e^{-\lam_1 t}+\rho_1^2,\hs  \,t\geq 0,
\ee
  \el
\noindent{{\bf Proof.}} The proof for  \eqref{e:3.1} is almost the same as in the
proof of \cite[Theorem 3.1]{Li2}. We omit the details and  focus our attention on the second estimate   \eqref{e:3.6}.

Observe  that $-\De u_k =\sum^k_{j=1}\mu_j a_{kj}(t)w_j$. Multiplying \eqref{e:2.3} by $\mu_j a_{kj}(t)$ and
summing the results for $j=1,\ldots,k$, we obtain that
\be\label{e:3.23}\begin{split}
\frac{1}{2}\frac{d}{dt}\| u_k\|_{V_1}^{2}+\|u_k\|_{V_2}^2&=-\( F_k(t,u_k,(u_k)_t),\,\De u_k\)\\[1ex]
 &\leq \frac{1}{2}\|u_k\|_{V_2}^{2}+ \frac{1}{2}| F_k(t,u_k,(u_k)_t)|_2^2\\[1ex]
 &\leq \frac{1}{2}\|u_k\|_{V_2}^{2}+ C_1\(| f(u_k)|_2^2+|g(u_k,u_k(t-\tau(t,(u_k)_t))|_2^2+| h_k|_2^2\),
 \end{split}
 \ee
where (and below) $C$ and  ${C_i}'s$  denote general positive  constants independent of $k$ and the initial data.

By (H0) we deduce that
 \bes
 | f(u_k)|_2^2\leq a_0^2\int_\W(|u_k|^p+1)^2dx\leq  C_2\(|u_k|_q^{2p}+1\),
  \ees
  and
   \bes\begin{split}
  |g(u_k,u_k(t-\tau(t,(u_k)_t))|_2^2&\leq b_0^2\int_\W\(\ba{ll} |u_k|^{\b}+|u_k(t-\tau(t,(u_k)_t)|^{{\b}}+1\ea\)^2dx\\
  &\leq C_3\(\|(u_k)_t\|_{\cL_{L^q(\W)}^\8}^{2\b}+1\).
 \end{split} \ees
Substituting  these estimates into \eqref{e:3.23} we arrive at
\be\label{e:3.24}
\frac{d}{dt}\|u_k\|_{V_1}^{2}+ \|u_k\|_{V_2}^{2}\leq C_4\(|u_k|_q^{2p}+\|(u_k)_t\|_{\cL_{L^q(\W)}^\8}^{2\b}+1\).
\ee
 Hence by \eqref{e:3.1} one concludes that there exist $M,\lam'>0$ such that
\bes\label{e:3.25}\begin{split}
\frac{d}{dt}\|u_k\|_{V_1}^{2}+ \|u_k\|_{V_2}^{2}\leq M\|\phi_k\|_{\cL_{L^q(\W)}^\8}^qe^{-\lam' t}+C_5,\Hs t>0.
 \end{split}
 \ees
By virtue of the Poinc\'{a}re inequality  we therefore have
\bes\label{e:3.26}\begin{split}
\frac{d}{dt}\|u_k\|_{V_1}^{2}&\leq - \|u_k\|_{V_2}^{2}+M\|\phi_k\|_{\cL_{L^q(\W)}^\8}^qe^{-\lam' t}+C_5\\[1ex]
&\leq - \mu_1\|u_k\|_{V_1}^{2}+M\|\phi_k\|_{\cL_{L^q(\W)}^\8}^qe^{-\lam_1 t}+C_5,\hs t>0,
 \end{split}
 \ees
where $\lam_1=\min\{\lam',\, \mu_1/2\}$. Thus  by the classical Gronwall lemma one deduces  that  there exist  $ M_1, \rho_1>0$ such   that
\be\label{e:3.27}\begin{split}
\|u_k\|_{V_1}^{2}&\leq \|\phi_k\|_{\cC_{ V_1}}^2e^{- \mu_1 t}+ \frac{M}{\mu_1-\lam_1}\|\phi_k\|_{\cL_{L^q(\W)}^\8}^{q}\(e^{- \lam_1 t}-e^{-\mu_1 t}\)+\rho_1^2\\[1ex]
&\leq \|\phi_k\|_{\cC_{ V_1}}^2 e^{- \mu_1 t}+M_1\|\phi_k\|_{\cL_{L^q(\W)}^\8}^{q}e^{- \lam_1 t}+\rho_1^2,\Hs t\geq 0.\\
\end{split}
 \ee
Hence  we  finish  the proof of the lemma.  \hfill $\Box$

\begin{remark}\label{r:1.2}
It can be assumed that  $M_0 > 1$. Then  it is trivial to see that \eqref{e:3.1} and \eqref{e:3.6} remain valid for $t\in [-r,0]$ since $u_k|_{[-r,0]}= \phi_k$.
\end{remark}


\br\label{r:1.3}
Let $T>0$. Integrating \eqref{e:3.24} between $0$ and $T$,  by \eqref{e:3.6} it yields
\be\label{e:3.7}
\begin{split}
 \int_0^T\|u_k\|_{V_2}^2\,ds&\leq C_T\(\|\phi_k\|_{\cC_{\bV_1}}^2+\|\phi_k\|_{\cL_{L^q(\W)}^\8}^{q}+1\).
\end{split}
\ee
\er

\subsection{Global existence}
For convenience in statement, we  write
\be\label{e:1.4}
\ol q=\frac{q-1}{\b_0}+1\hs \mb{for $q>1$}.\ee

\bt\label{t:3.1}
Assume $f$ and $g$ satisfy (H0)-(H1), and that $\tau$ satisfies (H2). Let $q_{c}<q<\8$, and  $ h\in L^{\8}(\mathbb{R};L^{\ol q}(\W))\cap C(\mathbb{R};L^{\ol q }(\W))$. Then for each $\phi\in \mathscr{X}_1^q$, system \eqref{e:1.1} has a global solutions $u=u(t;\phi)$ with
\be\label{e:1.22}
\left\{\ba{ll}  u\in C\([-r,\8); V_1\)+ L^\8\(-r,\8; V_1\)+ L^\8\(-r,\8;L^q(\W)\),\\[1.5ex] u\in L^2\(0,T;V_2\),\hs  u'\in L^2\(0,T;H\),\Hs \,T>0,\\[.5ex]\ea\right.\ee
Furthermore, there exist $B_i$, $\eta_i$ and  $\rho_i>0$ ($i=0,1$) such that
\be\label{e:q.3}
|u (t;\phi )|_{q}\leq B_0 \|\phi\|_{\cL_{L^q(\W)}^\8} e^{-\eta_0 t} +\rho_0,\Hs t\geq 0, \hs 1<q<\8,
\ee and
\be\label{e:q.5}
\| u (t;\phi )\|_{V_1} \leq B_1\(\|\phi\|_{\cC_{V_1}}+\|\phi\|_{\cL_{L^q(\W)}^\8}^{q/2}\) e^{-\eta_1 t}+\rho_1,\Hs  \,t\geq 0.
\ee
\et


\noindent{{\bf Proof.}} The proof for the existence part  is quite standard via Galerkin method. The key point is  to analyze the convergence of the  nonlinear term $\{F_k(t,u_k,(u_k)_t)\}$ for the Galerkin approximation solutions $u_k$ given above.

By virtue of \eqref{e:2.1} it can be assumed that
$$\|\phi_k\|_{\cC_{V_1}}\leq 2\|\phi\|_{\cC_{V_1}},\Hs k\geq 1.$$
We also recall that $\|\phi_k\|_{\mathcal{L}^{\8}_{L^q (\W)}}\leq 8 \|\phi\|_{\mathcal{L}^{\8}_{L^q (\W)}}$ for all $k$. Thus invoking  Lemma \ref{t:1.1} and Remark \ref{r:1.2} it is easy to deduce  that
there exists $C>0$ such that
\be\label{e:3.5}
\|u_k\|_{L^{\8}(-r,\8;V_1)}+\|u_k\|_{L^{\8}(-r,\8;L^q(\W))}\leq C.
\ee
As $q>q_{c}\geq \max\{2p,2\b\}$, by  (H0) it follows that
$$
|f(u_k)|_2+|g(u_k,u_k(t-\tau(t,(u_k)_t))|_2\leq C,\Hs t\geq 0,\,\,k\geq 1.$$

\vs

Let $T>0$ be fixed. Then   $\{u_k\}$ is uniformly bounded in $L^{\8}(-r,T;V_1)$.
Further  it can be shown that  $\{u'_k\}$ is uniformly bounded in $L^2(0,T;H)$.
Thanks to the classical  Lions-Aubin-Simon Lemma (see \cite[Corollary 4]{Simon}),  up to a subsequence,  $\{u_k\}$ converges in $C([0,T];X^{\gamma})$ to a function $u$ for   $0 \leq\gam< \frac{1}{2}$, where $X^\gam=D(A^{\gam})$ is the fractional power space. Consequently
\be\label{e:3.17a}u_k\ra u\mb{ in  $L^2(Q_T)$},\hs Q_T=\W\times(0,T).\ee It also follows that  there is a subsequence of $\{u_k\}$ (still denoted by $\{u_k\}$) such that
\be\label{e:3.17}
u_k(x,t) \rightarrow u(x,t),\Hs \mb{{a.e.}}\,\,(x,t) \in Q_T.
\ee

Using \eqref{e:3.17a} and hypothesis (H2) one easily verifies that
\be\label{e:3.17c}
\tau(t,(u_k)_t)\ra \tau(t,u_t)\hs\mb{as }k\ra\8
\ee
uniformly for $t\in[0,T]$.

Now by very standard argument it can be shown that, up to a subsequence,
$\{f(u_k(x,t))\}$ and $\{g(u_k,u_k(t-\tau(t,(u_k)_t))\}$ weakly converges to $f(u(x,t))$ and $g(u,u(t-\tau(t,u_t))$ in $L^2(0,T;H)$, respectively.


Passing to the limit in \eqref{e:2.3} one immediately concludes that
$u$ is  a weak solution  of system \eqref{e:1.1} on $[-r,T)$. Since $T$ is arbitrary, $u$ is actually a global weak solution of \eqref{e:1.1}.

It is easy to see that all the estimates obtained for the Galerkin approximations $u_k$ remain valid for $u$.

The relations $u\in L^2(0,T;V_2)$ and $u'\in L^2(0,T;H)$ in \eqref{e:1.22} are consequences of Remark \ref{r:1.3} and the boundedness of $\{u'_k\}$ in $L^2(0,T;H)$, respectively.

The continuity of $u$ in $V_1$ follows from  \cite[Chap. II, Theorem 3.3]{Tem} on abstract linear equations.
The proof of the theorem is complete.   \hfill $\Box$

\subsection{$L^{\infty}$-smoothing property and global $L^{\infty}$-estimate}
Based on these estimates of weak solutions $u$, we can obtain an instantaneous $L^\infty$-smoothing property of \eqref{e:1.1}.

\bt\label{t:1.2x} Assume that hypotheses {\em (H0)-(H2)} are fulfilled. Suppose
$\bh\in  L^{\8}(\mathbb{R};L^{\8}(\W))\cap C(\mathbb{R};L^{\8}(\W))$. Let $q_{c}^*<q<\8$ and $\hat p=\max\{p,\b\}$.
Then  there exist $B_2,\eta_2,\rho_2>0$ such that
\be\label{e:3.2}
|u(t;\phi)|_\8\leq B_2\|\phi\|_{\cL^\8_{L^q(\W)}}^{\hat p}t^{-\frac{d+1}{2q}}e^{-\eta_2 t }+\rho_2,\Hs t\geq0
\ee
for all solutions $u(t;\phi)$ of \eqref{e:1.1} with $\phi\in\sX_1^q$, where $d$ is the space dimension.
\et

\noindent{{\bf Proof.}} The argument is almost the same as in the proof of \cite[Theorem 3.6]{Li2}  by applying some fundamental $L^q - L^\infty$ estimates for linear equations (which can be found in Quittner and Souplet \cite[pp. 556]{QS}), we therefore omit the details.  \hfill $\Box$

\vs
 Repeating a similar  argument as in the proof of  \cite[Theorem 1.2]{Wang1} with almost no  modifications, we can obtain a global $L^{\8}$-estimate of the solution $u$  by letting $q\ra\8$ in \eqref{e:q.3}. Specifically, we have
\bt\label{t:1.3}
Assume $f$ and $g$ satisfy {\em (H0)-(H2)}, and let $ h\in  L^{\8}(\mathbb{R};L^{\8}(\W))\cap C(\mathbb{R};L^{\8}(\W))$. Then there exists $\rho_* >0$ such that for all $\phi\in \sX_1^{\8}$, we have
\be\label{e:x1}
|u(t;\phi)|_{\8}\leq \|\phi\|_{\cL^{\8}_{L^{\8}(\Omega)}}+\rho_*,\Hs t\geq 0.
\ee
\et

Fixing a $q>q_c^*$ and combining the above two estimates together, one immediately obtains a global $L^\8$-decay estimate for $u(t;\phi)$:

\bt\label{t:3.7} Assume that hypotheses {\em (H0)-(H2)} are fulfilled. Suppose
$\bh\in  L^{\8}(\mathbb{R};L^{\8}(\W))\cap C(\mathbb{R};L^{\8}(\W))$. Let $q_{c}^*<q<\8$ and $\hat p=\max\{p,\b\}$.
Then  there exist $B_3,\rho_\8>0$ such that
\be\label{e:3.2}
|u(t;\phi)|_\8\leq B_3\|\phi\|_{\cL^\8_{L^\8(\W)}}^{\hat p}e^{-\eta_2 t }+\rho_\8,\Hs t\geq0
\ee
for all solutions $u(t;\phi)$ of \eqref{e:1.1} with $\phi\in\sX_1^\8$.
\et

\Vs
\end{document}